# An approach for designing a surface pencil through a given geodesic curve

Gülnur SAFFAK ATALAY, Fatma GÜLER, Ergin BAYRAM*, Emin KASAP
Ondokuz Mayıs University, Faculty of Arts and Sciences, Mathematics Department
gulnur.saffak@omu.edu.tr, f.guler@omu.edu.tr ,erginbayram@yahoo.com,
kasape@omu.edu.tr

## ABSTRACT

Surfaces and curves play an important role in geometric design. In recent years, problem of finding a surface passing through a given curve have attracted much interest. In the present paper, we propose a new method to construct a surface interpolating a given curve as the geodesic curve of it. Also, we analyze the conditions when the resulting surface is a ruled surface. In addition, developablity along the common geodesic of the members of surface family are discussed. Finally, we illustrate this method by presenting some examples.

## 1. Introduction

A rotation minimizing adapted frame (RMF) {T,U,V} of a space curve contains the curve tangent T and the normal plane vectors U, V which show no instantaneous rotation about T. Because of their minimum twist RMFs are very interesting in computer graphics, including free-form deformation with curve constraints [1 - 6], sweep surface modeling [7 - 10], modeling of generalized cylinders and tree branches [11 - 15], visualization of streamlines and tubes [15 - 17], simulation of ropes and strings [18], and motion design and control [19].

There are infinitely many adapted frames on a given space curve [20]. One can produce other adapted frames from an existing one by controlling the orientation of the frame vectors U and V in the normal plane of the curve. In differential geometry the most familiar adapted frame is Frenet frame {T, N, B}, where T is the curve tangent, N is the principal normal vector and $B = T \times N$ is the binormal vector (see [21] for details). Beside of its fame, the Frenet frame is not a RMF and it is unsuitable for specifying the orientation of a rigid body along a given curve in applications such as motion planning, animation, geometric design, and robotics, since it incurs “unnecessary” rotation of the body [22]. Furthermore, Frenet frame is undefined if the curvature vanishes.

One of most significant curve on a surface is geodesic curve. Geodesics are important in the relativistic description of gravity. Einstein’s principle of equivalence tells us that geodesics represent the paths of freely falling particles in a given space. (Freely falling in this context means moving only under the influence of gravity, with no other forces involved). The geodesics principle states that the free trajectories are the geodesics of space. It plays a very important role in a geometric-relativity theory, since it means that the fundamental equation of dynamics is completely determined by the geometry of space, and therefore has not to be set as an independent equation. Moreover, in such a theory the action identifies (up to a constant) with the fundamental length invariant, so that the stationary action principle and the geodesics principle become identical. The concept of geodesic also finds its place in various industrial applications, such as tent manufacturing, cutting and painting path, fiberglass tape windings in pipe manufacturing, textile manufacturing [23–28]. In architecture, some special curves have nice properties in terms of structural functionality and manufacturing cost. One example is planar curves in vertical planes, whichcan be used as support elements. Another example is geodesic curves, [29]. Deng, B. , described methods to create patterns of special curves on surfaces, which find applications in design and realization of freeform architecture. He presented an evolution approach to generate a series of curves which are either geodesic or piecewise geodesic, starting from a given source curve on a surface. In [29], he investigated families of special curves (such as geodesics) on freeform surfaces, and propose computational tools to create such families. Also, he investigated patterns of special curves on surfaces, which find applications in design and realization of freeform architectural shapes (for details, see [29] ). Most people have heard the phrase; a straight line is the shortest distance between two points. But in differential geometry, they say this same thing in a different language. They say instead geodesics for the Euclidean metric are straight lines. A geodesic is a curve that represents the extremal value of a distance function in some space. In the Euclidean space, extremal means 'minimal',so geodesics are paths of minimal arc length. In general relativity, geodesics generalize the notion of "straight lines" to curved space time. This concept is based on the mathematical concept of a geodesic. Importantly, the world line of a particle free from all external force is a particular type of geodesic. In other words, a freely moving particle always moves along a geodesic. Geodesics are curves along which geodesic curvature vanishes. This is of course where the geodesic curvature has its name from.

In recent years, fundamental research has focused on the reverse problem or backward analysis: given a 3D curve, how can we characterize those surfaces that possess this curve as a special curve, rather than finding and classifying curves on analytical curved surfaces. The concept of family of surfaces having a given characteristic curve was first introduced by Wang et.al. [28] in Euclidean 3-space. Kasap et.al. [30] generalized the work of Wang by introducing new types of marching-scale functions, coefficients of the Frenet frame appearing in the parametric representation of surfaces. Also, surfaces with common geodesic in Minkowski 3-space have been the subject of many studies. In [31] Kasap and Akyıldız defined surfaces with a common geodesic in Minkowski 3-space and gave the sufficient conditions on marching-scale functions so that the given curve is a common geodesic on that surfaces. Şaffak and Kasap [32] studied family of surfaces with a common null geodesic**.** Lie et al. [33] derived the necessary and sufficient condition for a given curve to be the line of curvature on a surface. Bayram et al. [34] studied parametric surfaces which possess a given curve as a common asymptotic. However, they solved the problem using Frenet frame of the given curve.

In this paper, we obtain the necessary and sufficient condition for a given curve to be both isoparametric and geodesic on a parametric surface depending on the RMF. Furthermore, we show that there exists ruled surfaces possessing a given curve as a common geodesic curve and present a criteria for these ruled surfaces to be developable ones. We only study curves with an arc length parameter because such a study is easy to follow; if necessary, one can obtain similar results for arbitrarily parameterised regular curves.

## 2. Backgrounds

A parametric curve $r(s)$, $L_1 \leq s \leq L_2$, is a curve on a surface $P = P(s,t)$ in $\mathbb{R}^3$ that has a constant $s$ or *t*-parameter value. In this paper, $r'(s)$ denotes the derivative of *r* with respect to arc length parameter *s* and we assume that $r(s)$ is a regular curve, i.e. $r'(s) \neq 0$. For every point of $r(s)$, if $r''(s) \neq 0$, the set $\{T(s), N(s), B(s)\}$ is called the Frenet frame along $r(s)$, where $T(s) = r'(s)$, $N(s) = r''(s)/|r''(s)|$ and $B(s) = T(s) \times N(s)$ are the unit tangent, principal normal, and binormal vectors of the curve at the point $r(s)$, respectively. Derivative formulas of the Frenet frame is governed by the relations

$$\frac{d}{ds}\begin{pmatrix} T(s) \\ N(s) \\ B(s) \end{pmatrix} = \begin{pmatrix} 0 & \kappa(s) & 0 \\ -\kappa(s) & 0 & \tau(s) \\ 0 & -\tau(s) & 0 \end{pmatrix}\begin{pmatrix} T(s) \\ N(s) \\ B(s) \end{pmatrix} \tag{2.1}$$

where $\kappa(s) = \|r''(s)\|$ and $\tau(s) = \dfrac{(r'(s), r''(s), r'''(s))}{\|r'(s) \times r''(s)\|}$ are called the curvature and torsion of the curve $r(s)$, respectively [35].

Another useful frame along a curve is rotation minimizing frame. They are useful in animation, motion planning, swept surface constructions and related applications where the Frenet frame may prove unsuitable or undefined. A frame $\{T(s), U(s), V(s)\}$ among the frames on the curve $r(s)$ is called *rotation minimizing* if it is the frame of minimum twist around the tangent vector *T*. $\{T(s), U(s), V(s)\}$ is an RMF if

$$\begin{cases} U'(s) = -(U(s) \bullet r''(s))r'(s), \\ V'(s) = -(V(s) \bullet r''(s))r'(s), \end{cases}$$

where "•" denotes the standard inner product in $\mathbb{R}^3$ [36]. Observe that such a pair *U* and *V* is not unique; there exist a one parameter family of RMF's corresponding to different sets of initial values of *U* and *V*. According to Bishop [20], a frame is an RMF if and only if each of $U'(s)$ and $V'(s)$ is parallel to $T(s)$. Equivalently,

$$\begin{cases} U'(s) \bullet V(s) \equiv 0, \\ V'(s) \bullet U(s) \equiv 0 \end{cases} \tag{2.2}$$

is the necessary and sufficient condition for the frame to be rotation minimizing [37].
There is a relation between the Frenet frame (if the Frenet frame is defined) and RMF, that is, *U* and *V* are the rotation of *N* and *B* of the curve $r(s)$ in the normal plane. Then,

$$\begin{bmatrix} U \\ V \end{bmatrix} = \begin{bmatrix} \cos\theta & \sin\theta \\ -\sin\theta & \cos\theta \end{bmatrix}\begin{bmatrix} N \\ B \end{bmatrix}, \tag{2.3}$$

where $\theta = \theta(s)$ is the angle between the vectors *N* and *U* (see Fig. 2), [38].

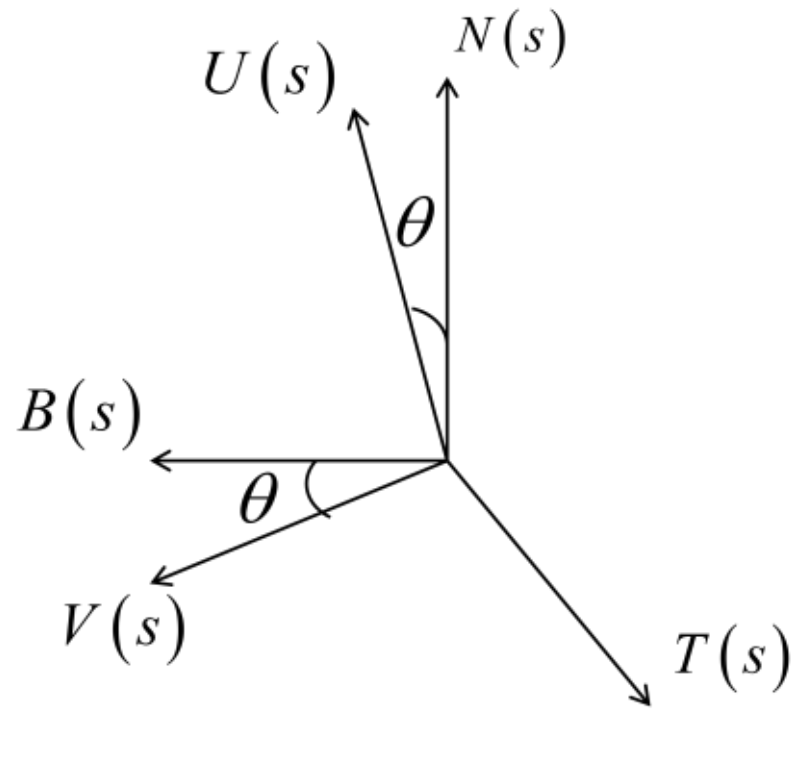

**Fig. 2** The Frenet frame (*T(s), N(s), B(s)*) and the vectors *U(s), V(s).*

Eqn. (2.3) implies that $\{T(s), U(s), V(s)\}$ is an RMF if it satisfies the following relations

$$U' = -\kappa\cos\theta T \quad, V' = \kappa\sin\theta T, \;\; \theta' = -\tau . \tag{2.4}$$

Note that $\{T(s), U(s), V(s)\}$ is defined along the curve $r(s)$ even if the curvature vanishes where the Frenet frame is undefined.

### 3. Surface pencil with an geodesic curve

Suppose we are given a 3-dimensional parametric curve $r(s),\ L_1 \le s \le L_2,$ in which $s$ is the arc length (regular and $\|r'(\mathrm{s})\| = 1,\ L_1 \le s \le L_2$).

Surface pencil that interpolates $r(s)$ as a common curve is given in the parametric form as

$$P(s,t) = r(s) + a(s,t)T(s) + b(s,t)U(s) + c(s,t)V(s),\ L_1 \le s \le L_2,\ T_1 \le t \le T_2, \quad (3.1)$$

where $a(s,t), b(s,t)$ and $\mathrm{c}(s,t)$ are $C^1$ functions. The values of the *marching-scale functions* $a(s,t), b(s,t)$ and $\mathrm{c}(s,t)$ indicate, respectively, the extension-like, flexion-like, and retortion-like effects caused by the point unit through time *t*, starting from $r(s)$.

**Remark 3.1 :** Observe that choosing different marching-scale functions yields different surfaces possessing $r(s)$ as a common curve.

Our goal is to find the necessary and sufficient conditions for which the curve $r(s)$ is isoparametric and geodesic on the surface $P(s,t)$. Firstly, as $r(s)$ is an isoparametric curve on the surface $P(s,t)$, there exists a parameter $t_0 \in [T_1, T_2]$ such that

$$a(s,t_0) = b(s,t_0) = c(s,t_0) \equiv 0,\ L_1 \le s \le L_2,\ T_1 \le t_0 \le T_2. \quad (3.2)$$

Secondly the curve $r(s)$ is an geodesic curve on the surface $P(s,t)$ there exist a parameter $t_0 \in [T_1, T_2]$ such that

According to the geodesic theory [39], geodesic curvature $k_g = \det(r', r'', n)$ vanishes along geodesics. Thus, if we get:

$$n(s,t_0) \,//\, \mathrm{N(s)}. \quad (3.3)$$

where $n(s,t_0)$ is the surface normal along the curve $r(s)$ and N is a normal vector of $r(s)$

The normal vector of $P = P(s,t)$ can be written as

$$n(s,t) = \frac{\partial P(s,t)}{\partial s} \times \frac{\partial P(s,t)}{\partial t}.$$

From Eqns. (2.1) and (2.3), the normal vector can be expressed as

$$\begin{aligned} n(s,t) = {} & \left( \frac{\partial c(s,t)}{\partial t} \left( a(s,t)\kappa(s)\cos\theta(s) + \frac{\partial b(s,t)}{\partial s} \right) - \frac{\partial b(s,t)}{\partial t} \left( \frac{\partial c(s,t)}{\partial s} - a(s,t)\kappa(s)\sin\theta(s) \right) \right) T(s) \\ & + \left( \frac{\partial a(s,t)}{\partial t} \left( \frac{\partial c(s,t)}{\partial s} - a(s,t)\kappa(s)\sin\theta(s) \right) - \frac{\partial c(s,t)}{\partial t} \left( 1 + \frac{\partial a(s,t)}{\partial s} - b(s,t)\kappa(s)\cos\theta(s) + c(s,t)\kappa(s)\sin\theta(s) \right) \right) U(s) \\ & + \left( \frac{\partial b(s,t)}{\partial t} \left( 1 + \frac{\partial a(s,t)}{\partial s} - b(s,t)\kappa(s)\cos\theta(s) + c(s,t)\kappa(s)\sin\theta(s) \right) - \frac{\partial a(s,t)}{\partial t} \left( a(s,t)\kappa(s)\cos\theta(s) + \frac{\partial b(s,t)}{\partial s} \right) \right) V(s) \end{aligned}$$

Thus,

$$n(s,t_0)=\phi_1(s,t_0)T(s)+\phi_2(s,t_0)U(s)+\phi_3(s,t_0)V(s) \quad (3.4)$$

where

$$\phi_1(s,t_0)=\left(\frac{\partial c}{\partial t}(s,t_0)\left(a(s,t_0)\kappa(s)\cos\theta(s)+\frac{\partial b}{\partial s}(s,t_0)\right)-\frac{\partial b}{\partial t}(s,t_0)\left(\frac{\partial c}{\partial s}(s,t_0)-a(s,t_0)\kappa(s)\sin\theta(s)\right)\right),$$

$$\phi_2(s,t_0)=\left(\frac{\partial a}{\partial t}(s,t_0)\left(\frac{\partial c}{\partial s}(s,t_0)-a(s,t_0)\kappa(s)\sin\theta(s)\right)-\frac{\partial c}{\partial t}(s,t_0)\left(1+\frac{\partial a}{\partial s}(s,t_0)-b(s,t_0)\kappa(s)\cos\theta(s)+c(s,t_0)\kappa(s)\sin\theta(s)\right)\right),$$

$$\phi_3(s,t_0)=\left(\frac{\partial b}{\partial t}(s,t_0)\left(1+\frac{\partial a}{\partial s}(s,t_0)-b(s,t_0)\kappa(s)\cos\theta(s)+c(s,t_0)\kappa(s)\sin\theta(s)\right)-\frac{\partial a}{\partial t}(s,t_0)\left(a(s,t_0)\kappa(s)\cos\theta(s)+\frac{\partial b}{\partial s}(s,t_0)\right)\right).$$

**Remark 3.2 :** Because,

$$\begin{cases} a(s,t_0)=b(s,t_0)=c(s,t_0)\equiv 0, \\ t_0\in[T_1,T_2],\ L_1\le s\le L_2, \end{cases}$$

along the curve $r(s)$, by the definition of partial differentiation we have

$$\begin{cases} \frac{\partial a}{\partial s}(s,t_0)=\frac{\partial b}{\partial s}(s,t_0)=\frac{\partial c}{\partial s}(s,t_0)\equiv 0, \\ t_0\in[T_1,T_2],\ L_1\le s\le L_2. \end{cases}$$

According to remark above, we should have

$$\begin{cases} \phi_1(s,t_0)\equiv 0, \\ \phi_2(s,t_0)=-\frac{\partial c}{\partial t}(s,t_0), \\ \phi_3(s,t_0)=\frac{\partial b}{\partial t}(s,t_0). \end{cases} \quad (3.5)$$

Thus, from (3.4) we obtain

$$n(s,t_0)=\phi_2(s,t_0)U(s)+\phi_3(s,t_0)V(s)$$

from Eqn. (2.3), we get

$$n(s,t_0)=\left(\cos\theta(s)\phi_2(s,t_0)-\sin\theta(s)\phi_3(s,t_0)\right)N(s)+\left(\sin\theta(s)\phi_2(s,t_0)+\cos\theta(s)\phi_3(s,t_0)\right)B(s).$$

from Eqn. (3.3), we know that $r(s)$ is a geodesic curve if and only if

$$\cos\theta(s)\phi_2(s,t_0)-\sin\theta(s)\phi_3(s,t_0)\neq 0 \quad (3.6)$$

and

$$\sin\theta(s)\phi_2(s,t_0)+\cos\theta(s)\phi_3(s,t_0)=0 \quad (3.7)$$

From (3.6) and (3.7), we obtain

$$-\frac{1}{\sin\theta(s)}\phi_3(s,t_0)\neq 0$$

there using (3.5),we have

$\frac{\partial b}{\partial t}(s,t_0)\neq 0$ and $\sin\theta\frac{\partial c}{\partial t}(s,t_0)=\cos\theta\frac{\partial b}{\partial t}(s,t_0)$.

Thus we have following theorem:

**Theorem 3.3 :** The necessary and sufficient condition for the curve $r(s)$ to be both isoparametric and geodesic on the surface $P(s,t)$ is

$$\begin{cases} a(s,t_0)=b(s,t_0)=c(s,t_0)\equiv 0, \\ \sin\theta\frac{\partial c}{\partial t}(s,t_0)=\cos\theta\frac{\partial b}{\partial t}(s,t_0), \frac{\partial b}{\partial t}(s,t_0)\neq 0 \\ \theta'(s)=-\tau(s). \end{cases} \tag{3.8}$$

**Corollary 3.4 :** The sufficient condition for the curve $r(s)$ to be both isoparametric and geodesic on the surface $P(s,t)$ is

$$b(s,t)=f(s,t)\sin\theta,\ c(s,t)=f(s,t)\cos\theta,\ f(s,t_0)\equiv 0,\ \theta'(s)=-\tau(s). \tag{3.9}$$

**4. Ruled surface pencil with a common asymptotic curve**

**Theorem 4.1 :** Given an arc-length curve $r(s)$, there exists a ruled surface $P(s,t)$ possessing $r(s)$ as a common geodesic curve.

**Proof :** Choosing marching-scale functions as

$$a(s,t)=(t-t_0)g(s),\ b(s,t)=(t-t_0)\sin\theta,\ c(s,t)=(t-t_0)\cos\theta \tag{3.10}$$

and $\theta'(s)=-\tau(s)$ Eqn. (3.1) takes the following form of a ruled surface

$$P(s,t)=r(s)+(t-t_0)\left[g(s)T(s)+\sin\theta(s)U(s)+\cos\theta(s)V(s)\right], \tag{3.11}$$

which satisfies Eqn. (3.8) interpolating $r(s)$ as a common geodesic curve.

**Corollary 4.2 :** Ruled surface (3.11) is developable if and only $\tau(s)=g(s)\kappa(s)$ is to be.

**Proof :** $P(s,t)=r(s)+(t-t_0)\left[g(s)T(s)+\sin\theta(s)U(s)+\cos\theta(s)V(s)\right]$ is developable if and only if $\left(r',d,d'\right)=0,$ where $d(s)=g(s)T(s)+\sin\theta(s)U(s)+\cos\theta(s)V(s)$. Using Eqns. (2.1) – (2.4) gives

$$\begin{aligned} d' &= g'T+gT'+\theta'\cos\theta U+\sin\theta U'-\theta'\sin\theta V+\cos\theta V' \\ &= g'T+g\kappa\left(\cos\theta U-\sin\theta V\right)-\tau\cos\theta U-\kappa\sin\theta\cos\theta T+\tau sin\theta V+\kappa\sin\theta\cos\theta T \\ &= g'T+\left(g\kappa\cos\theta-\tau\cos\theta\right)U+\left(-g\kappa\sin\theta+\tau\sin\theta\right)V. \end{aligned} \tag{3.12}$$

Employing Eqn. (3.12) in the determinant we get $\tau=g\kappa$, which complete proof.

**5. Examples of generating surfaces with a common geodesic curve**

**Example 5.1.** Let $r(s)=\left(\frac{3}{5}\sin(s),\frac{3}{5}\cos(s),\frac{4}{5}s\right)$ be a unit speed curve. Then it is easy to show that

$$\begin{cases} T(s) = \left(\dfrac{3}{5}\cos(s), -\dfrac{3}{5}\sin(s), \dfrac{4}{5}\right) \\ \kappa(s) = \dfrac{3}{5} \quad , \quad \tau(s) = \dfrac{-4}{5} \end{cases}$$

If we choose $\theta = \frac{4}{5}s$ +c and c=0

$$\begin{cases} U' = \left( -\dfrac{9}{25}\cos(s)\cos(\dfrac{4s}{5}), \dfrac{9}{25}\sin(s)\cos(\dfrac{4s}{5}), -\dfrac{12}{25}\cos(\dfrac{4s}{5}) \right) \\ V' = \left( \dfrac{9}{25}\cos(s)\sin(\dfrac{4s}{5}), -\dfrac{9}{25}\sin(s)\sin(\dfrac{4s}{5}), \dfrac{12}{25}\sin(\dfrac{4s}{5}) \right). \end{cases},$$

then Eqn. (2.5) is satisfied. By integration, we obtain

$$\begin{cases} U = \left( -\dfrac{1}{10}\sin(\dfrac{9s}{5}) - \dfrac{9}{10}\sin(\dfrac{s}{5}),\ -\dfrac{1}{10}\cos(\dfrac{9s}{5}) - \dfrac{9}{10}\cos(\dfrac{s}{5}),\ -\dfrac{3}{5}\sin(\dfrac{4s}{5}) \right), \\ V = \left( -\dfrac{1}{10}\cos(\dfrac{9s}{5}) + \dfrac{9}{10}\cos(\dfrac{s}{5}),\ -\dfrac{9}{10}\sin(\dfrac{s}{5}) + \dfrac{1}{10}\sin(\dfrac{9s}{5}),\ -\dfrac{3}{5}\cos(\dfrac{4s}{5}) \right) \end{cases}$$

Now, $\{T(s), U(s), V(s)\}$ is an RMF since is satisfies Eqn. (2.2). If we take $a(s,t) \equiv 0,\ b(s,t) = \sin(\frac{4s}{5})(\sin t - 1),\ c(s,t) = \cos(\frac{4s}{5})\cos t$ and $t_0 = \frac{\pi}{2}$, then Eqn. (3.5) is satisfied. Thus, we obtain a member of the surface pencil with a common geodesic curve $r(s)$ as

$$P_1(s,t) = \begin{pmatrix} \dfrac{3}{5}\sin(s) + (\sin(t)-1)\sin\left(\dfrac{4s}{5}\right)\left(-\dfrac{1}{10}\sin\left(\dfrac{9s}{5}\right) - \dfrac{9}{10}\sin\left(\dfrac{s}{5}\right)\right) + \cos\left(\dfrac{4s}{5}\right)\cos(t)\left(-\dfrac{1}{10}\cos\left(\dfrac{9s}{5}\right) + \dfrac{9}{10}\cos\left(\dfrac{s}{5}\right)\right), \\ \dfrac{3}{5}\cos(s) + (\sin(t)-1)\sin\left(\dfrac{4s}{5}\right)\left(-\dfrac{1}{10}\cos\left(\dfrac{9s}{5}\right) - \dfrac{9}{10}\cos\left(\dfrac{s}{5}\right)\right) + \cos\left(\dfrac{4s}{5}\right)\cos(t)\left(-\dfrac{9}{10}\sin\left(\dfrac{s}{5}\right) + \dfrac{1}{10}\sin\left(\dfrac{9s}{5}\right)\right), \\ \dfrac{4}{5}s - \dfrac{3}{5}\left(\sin^2\left(\dfrac{4s}{5}\right)(\sin(t)-1) + \cos^2\left(\dfrac{4s}{5}\right)\cos(t)\right) \end{pmatrix}$$

where $0 \le s \le 2\pi,\ 0 \le t \le 2\pi$ (Fig. 1).

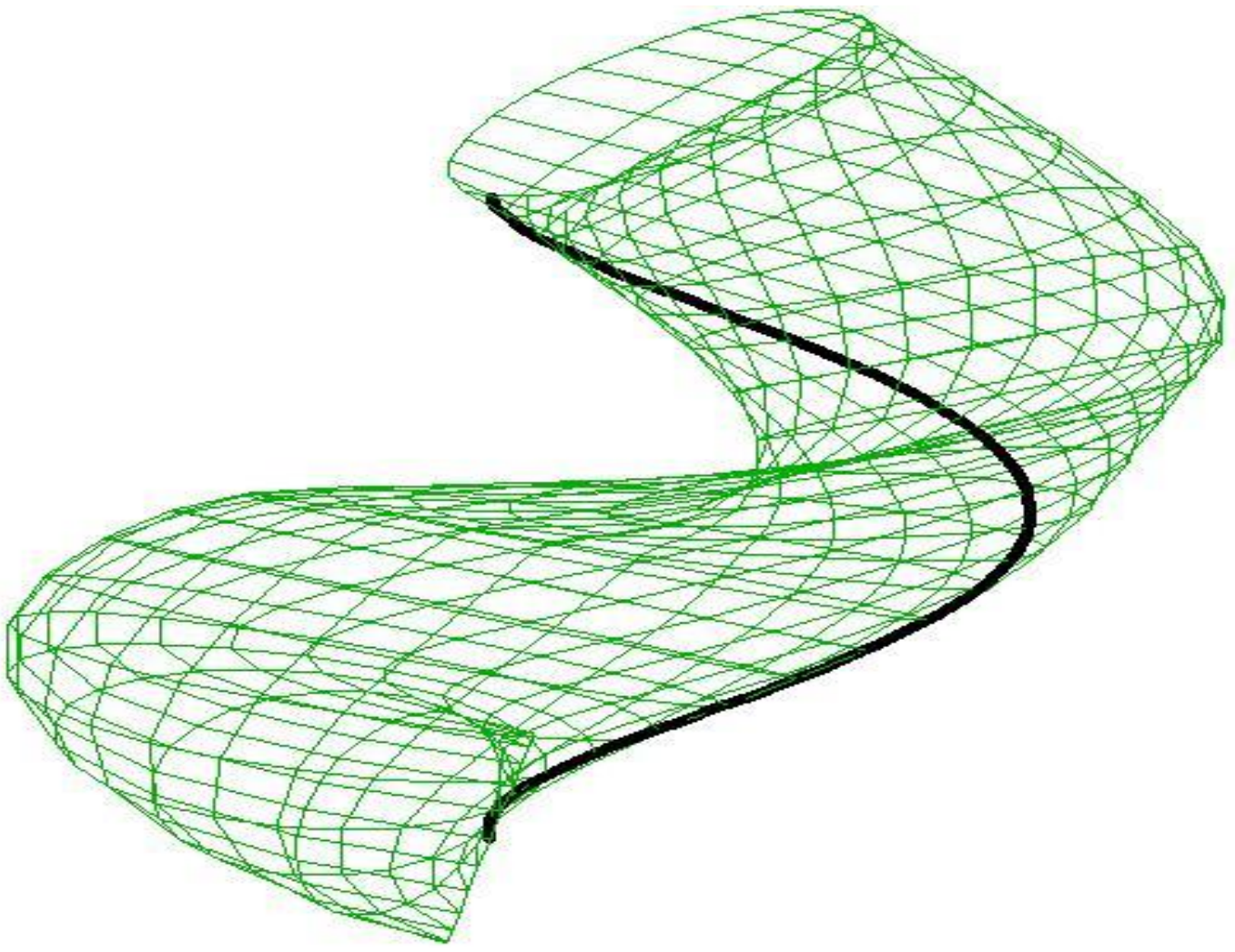

Fig. 1. $P_1(s,t)$ as a member of surface pencil and its geodesic curve.

In Eqn. (3.8), if we take $g(s)=\dfrac{\tau(s)}{\kappa(s)}=-\dfrac{4}{3}$, then we obtain the following developable ruled surface with a common geodesic curve $r(s)$ as

$$P_2(s,t)=\begin{pmatrix} \frac{3}{5}\sin(s)+t\left(-\frac{4}{5}\cos(s)+\sin\left(\frac{4s}{5}\right)\left(-\frac{1}{10}\sin\left(\frac{9s}{5}\right)-\frac{9}{10}\sin\left(\frac{s}{5}\right)\right)+\cos\left(\frac{4s}{5}\right)\left(-\frac{1}{10}\cos\left(\frac{9s}{5}\right)+\frac{9}{10}\cos\left(\frac{s}{5}\right)\right)\right), \\ \frac{3}{5}\cos(s)+t\left(\frac{4}{5}\sin(s)+\sin\left(\frac{4s}{5}\right)\left(-\frac{1}{10}\cos\left(\frac{9s}{5}\right)-\frac{9}{10}\cos\left(\frac{s}{5}\right)\right)+\cos\left(\frac{4s}{5}\right)\left(-\frac{9}{10}\sin\left(\frac{s}{5}\right)+\frac{1}{10}\sin\left(\frac{9s}{5}\right)\right)\right), \\ \frac{4}{5}s-\frac{5}{3}t \end{pmatrix}$$

where $-2\pi\le s\le 2\pi,\ \ 0\le t\le 2\pi$ (Fig. 2).

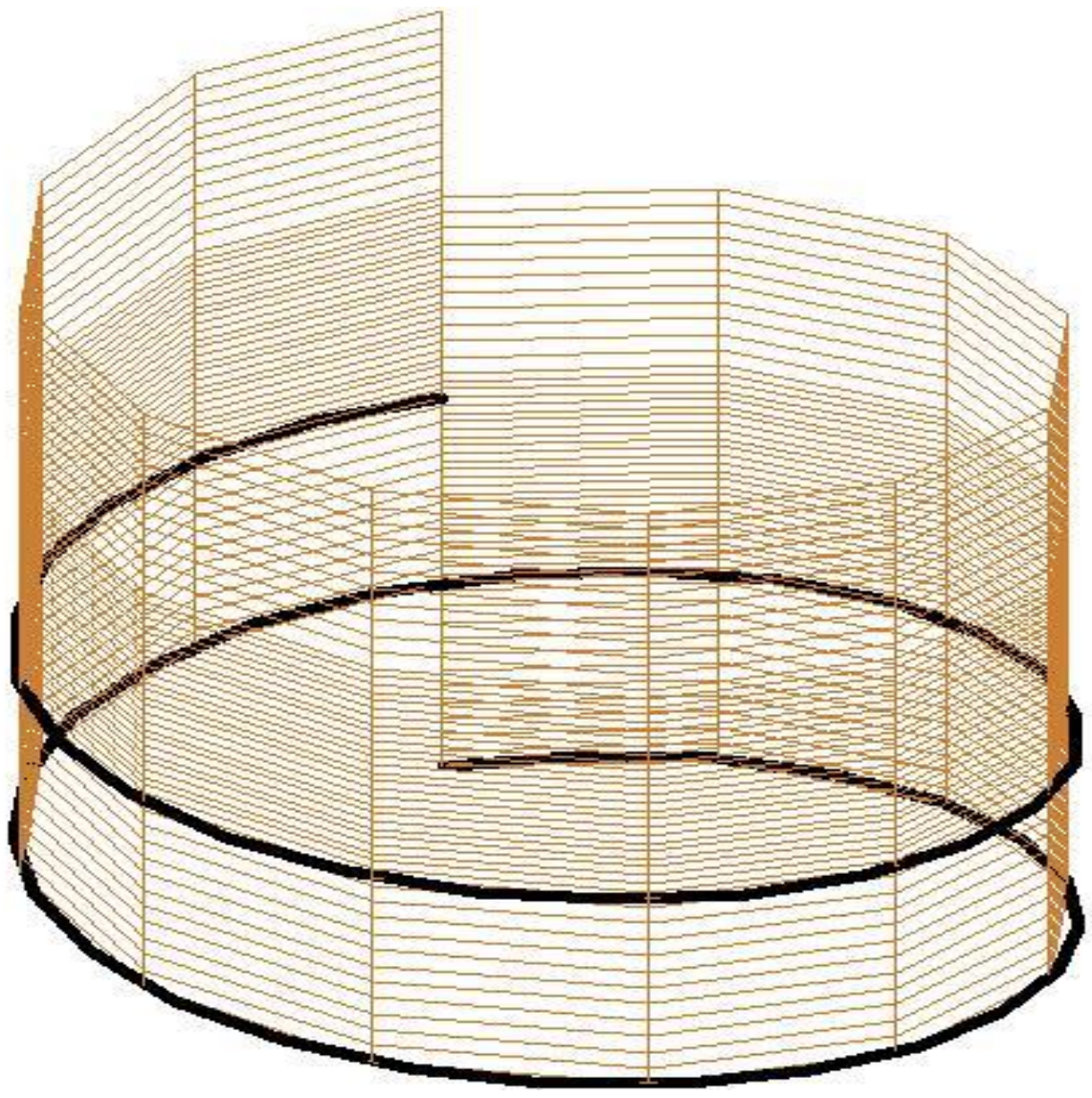

Fig. 2. $P_2(s,t)$ as a member of the developable ruled surface pencil and its geodesic curve.

**Example 5.2.** Let $r(s)=(\cos(s),\sin(s),0)$ be a unit speed curve. It is obvious that

$$T(s)=(-\sin(s),\cos(s),0),\ \kappa(s)=1,\ \tau(s)=0.$$

If we take $U(s)=(-\cos(s),-\sin(s),0)$ and $V(s)=(0,0,1)$, then Eqn. (2.2) is satisfied and $\{T(s),U(s),V(s)\}$ is an RMF. By choosing marching-scale functions as $a(s,t)\equiv 0$, $b(s,t)=1-\cos t$, $c(s,t)=\sin t$ and $t_0=\theta=0$, then Eqn. (3.5) is satisfied. Thus, we immediately obtain a member of the surface pencil with a common geodesic curve $r(s)$ as

$$P_3(s,t)=(\cos(s)\cos(t),\sin(s)\cos(t),\sin(t))$$

where $0\le s\le 2\pi,\ 0\le t\le 2\pi$ (Fig. 3).

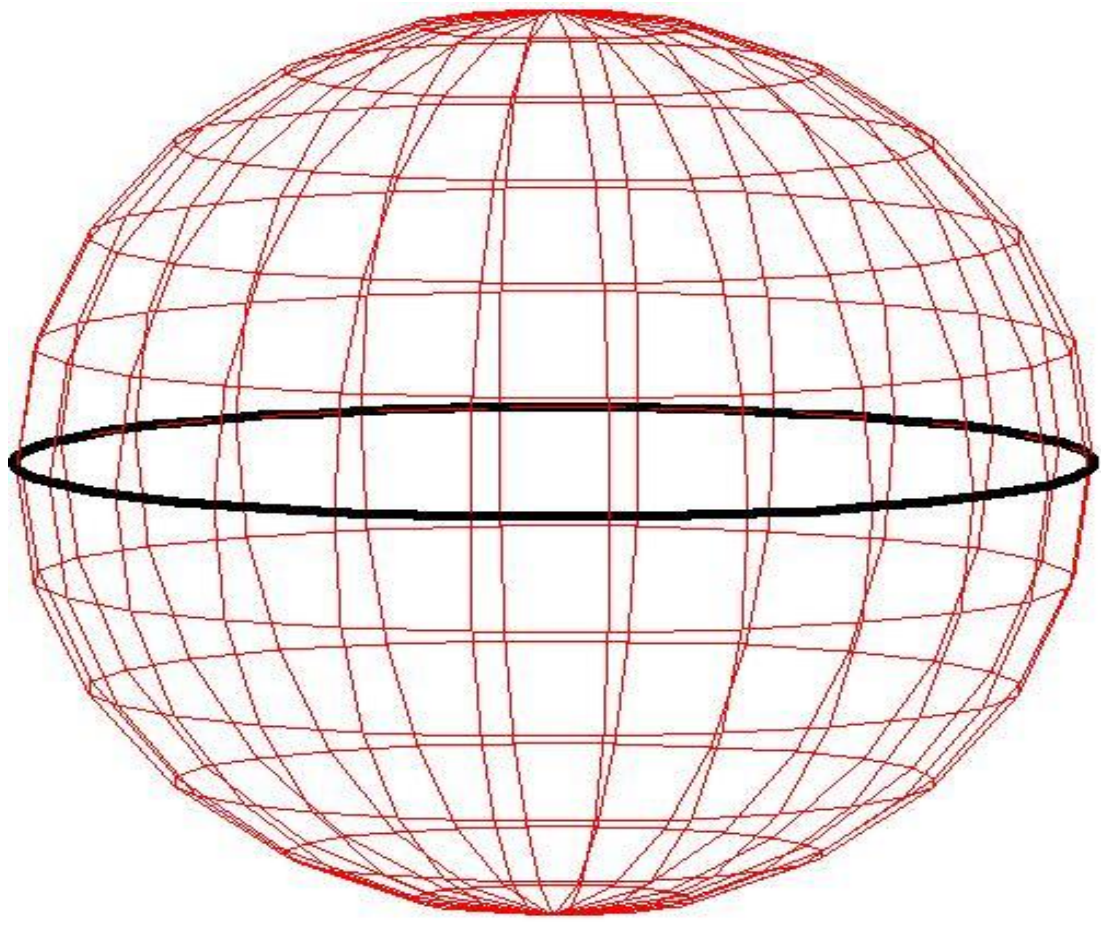

Fig. 3. $P_3(s,t)$ as a member of surface pencil and its geodesic curve.

For the same curve let us find a ruled surface. In Eqn. (3.8), if we take $g(s)=\dfrac{\tau(s)}{\kappa(s)}=0$, then we obtain the following developable ruled surface with a common geodesic curve $r(s)$ as

$$P_4(s,t)=(\cos(s),\sin(s),t)$$

where $-\pi\le s\le\pi$, $-1\le t\le 1$ (Fig. 4).

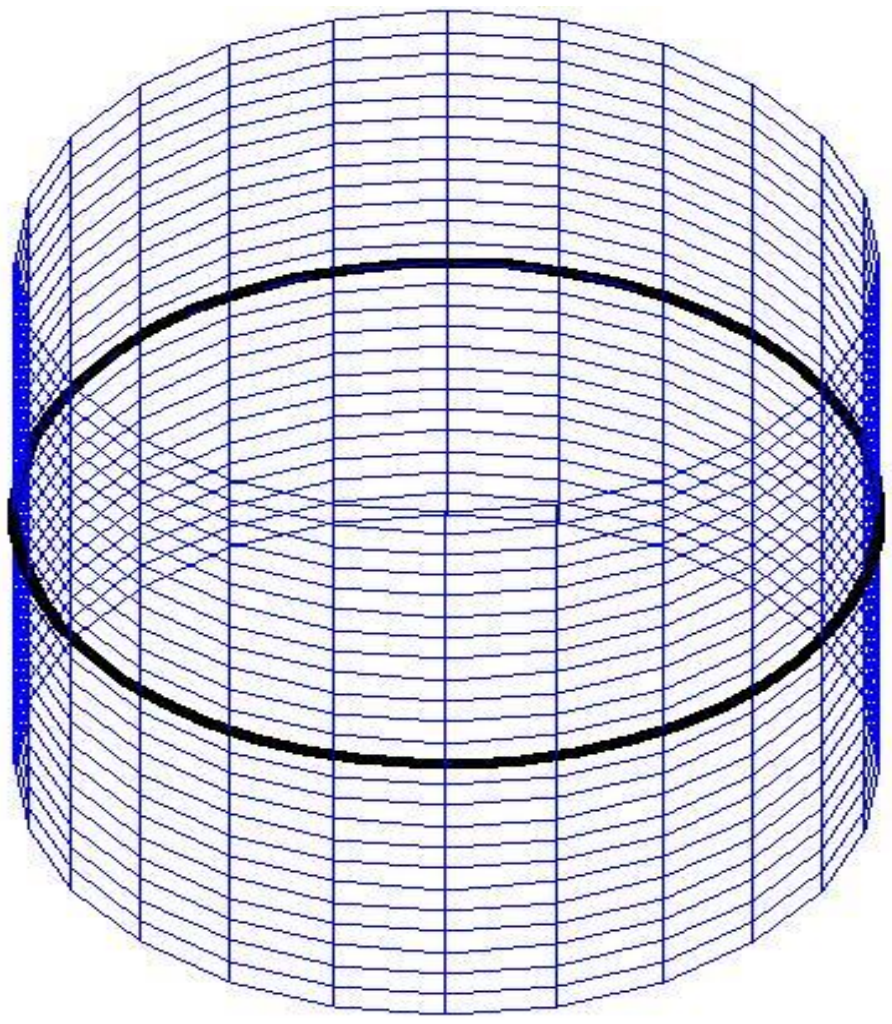

Fig. 4. $P_4(s,t)$ as a member of the developable ruled surface pencil and its geodesic curve.

If we take $g(s)$ = 0 and $t_0 = \theta = \frac{\pi}{3}$ then, we obtain the developable ruled surface with a common geodesic curve $r(s)$ as

$$P_5(s,t) = \left( \cos(s)\left(1 - \frac{\sqrt{3}}{2}\left(t - \frac{\pi}{3}\right)\right), \sin(s)\left(1 - \frac{-2+\sqrt{3}}{4}\left(t - \frac{\pi}{3}\right)\right), \left(\frac{2\sqrt{3}+1}{4}\left(t - \frac{\pi}{3}\right)\right)\right)$$

where $0 \le s \le 2\pi,\ 0 \le t \le 2\pi$ (Fig. 5).

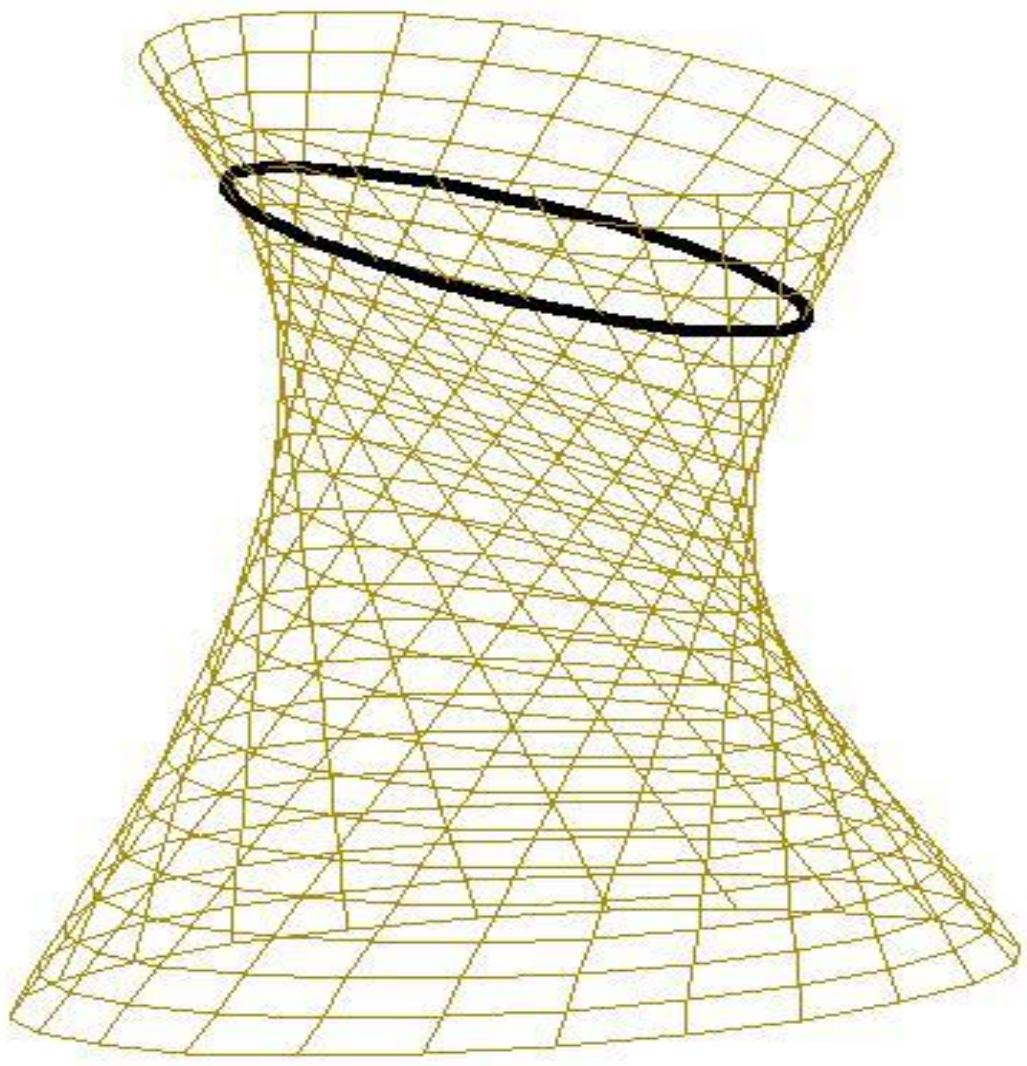

Fig. 5. $P_5(s,t)$ as a member of the developable ruled surface pencil and its geodesic curve.

## References

[1] Bechmann D., Gerber D. Arbitrary shaped deformation with dogme, Visual Comput. 19, 2–3, 2003, pp. 175-86.

[2] Peng Q, Jin X, Feng J. Arc-length-based axial deformation and length preserving deformation. In Proceedings of Computer Animation 1997; 86-92.

[3] Lazarus F, Coquillart S, Jancène P. Interactive axial deformations, In Modeling in Computer Graphics: Springer Verlag 1993;241-54.

[4] Lazarus F., Verroust A. Feature-based shape transformation for polyhedral objects. In Proceedings of the 5th Eurographics Workshop on Animation and Simulation 1994; 1-14.

[5] Lazarus F, Coquillart S, Jancène P. Axial deformation: an intuitive technique. Comput.-Aid Des 1994; 26, 8:607-13.

[6] Llamas I, Powell A, Rossignac J, Shaw CD. Bender: A virtual ribbon for deforming 3d shapes in biomedical and styling applications. In Proceedings of Symposium on Solid and Physical Modeling 2005; 89-99.

[7] Bloomenthal M, Riesenfeld RF. Approximation of sweep surfaces by tensor product NURBS. In SPIE Proceedings Curves and Surfaces in Computer Vision and Graphics II 1991; Vol 1610:132-54.

[8] Pottmann H, Wagner M. Contributions to motion based surface design. Int J Shape Model 1998; 4, 3&4:183-96.

[9] Siltanen P, Woodward C. Normal orientation methods for 3D offset curves, sweep surfaces, skinning. In Proceedings of Eurographics 1992; 449-57.

[10] Wang W, Joe B. Robust computation of rotation minimizing frame for sweep surface modeling. Comput.-Aid Des 1997; 29:379-91.

[11] Shani U, Ballard DH. Splines as embeddings for generalized cylinders. Comput Vision Graph Image Proces 1984; 27: 129-56.

[12] Bloomenthal J. Modeling the mighty maple. In Proceedings of SIGGRAPH 1985:305-11.

[13] Bronsvoort WF, Klok F. Ray tracing generalized cylinders. ACM Trans Graph 1985; 4, 4: 291–302.

[14] Semwal SK, Hallauer J. Biomedical modeling: implementing line-of-action algorithm for human muscles and bones using generalized cylinders. Comput Graph 1994:18, 1: 105-12.

[15] Banks DC, Singer BA. A predictor-corrector technique for visualizing unsteady flows. IEEE Trans on Visualiz Comput Graph 1995; 1, 2: 151-63.

[16] Hanson AJ, Ma H. A quaternion approach to streamline visualization. IEEE Trans Visualiz Comput Graph 1995; 1, 2: 164-74.

[17] Hanson A. Constrained optimal framing of curves and surfaces using quaternion gauss map. In Proceedings of Visulization 1998; 375-82.

[18] Barzel R. Faking dynamics of ropes and springs. IEEE Comput Graph Appl 1997; 17, 3: 31-39.

[19] Jüttler B. Rational approximation of rotation minimizing frames using Pythagorean-hodograph cubics. J Geom Graph 1999; 3: 141-59.

[20] Bishop RL. There is more than one way to frame a curve. Am Math Mon 1975; 82: 246-51.

[21] O'Neill B. Elementary Differential Geometry. New York: Academic Press Inc; 1966.

[22] Farouki RT, Sakkalis T. Rational rotation-minimizing frames on polynomial space curves of arbitrary degree. J Symbolic Comput 2010; 45: 844-56.

[23] Brond R, Jeulin D, Gateau P, Jarrin J, Serpe G. Estimation of the transport properties ofpolymer composites by geodesic propagation. J Microsc 1994;176:167–77.

[24] Bryson S. Virtual spacetime: an environment for the visualization of curved spacetimesvia geodesic flows. Technical Report, NASA NAS, Number RNR-92-009; March 1992.

[25] Grundig L, Ekert L, Moncrieff E. Geodesic and semi-geodesic line algorithms for cuttingpattern generation of architectural textile structures. In:Lan TT, editor. Proceedings of theAsia-Pacific Conference on Shell and Spatial Structures, Beijing. 1996.

[26] Haw RJ. An application of geodesic curves to sail design. Comput Graphics Forum1985;4(2):137–9.

[27] Haw RJ, Munchmeyer RC. Geodesic curves on patched polynomial surfaces. ComputGraphics Forum 1983;2(4):225–32.

[28] Wang G. J, Tang K, Tai C. L. Parametric representation of a surface pencil with a common spatial geodesic. Comput. Aided Des. 36 (5)(2004) 447-459.

[29] Deng, B. , 2011. Special Curve Patterns for Freeform Architecture Ph.D. thesis, Eingereicht an der Technischen Universitat Wien,Fakultat für Mathematik und Geoinformation von.

[30] Kasap E, Akyıldız FT, Orbay K. A generalization of surfaces family with common spatial geodesic. Appl Math Comput 2008; 201: 781–9.

[31] Kasap E, Akyildiz F.T. Surfaces with common geodesic in Minkowski 3-space. Applied Mathematics and Computation, 177 (2006) 260-270.

[32] Şaffak G, Kasap E. Family of surface with a common null geodesic. International Journalof Physical Sciences Vol. 4(8), pp. 428-433, August, 2009.

[33] Bayram E, Güler F, Kasap E. Parametric representation of a surface pencil with a common asymptotic curve. Comput Aided Des 2012; 44: 637-643.

[34] Li CY, Wang RH, Zhu CG. Parametric representation of a surface pencil with a common line of curvature. Comput Aided Des 2011;43(9):1110–7.

[35] do Carmo MP. Differential geometry of curves and surfaces. Englewood Cliffs (New Jersey): Prentice Hall, Inc.; 1976.

[36] Klok F. Two moving coordinate frames along a 3D trajectory. Comput Aided Geom Design 1986; 3:217-229.

[37] Han CY. Nonexistence of rational rotation-minimizing frames on cubic curves. Comput Aided Geom Design 2008; 25:298-304.

[38] Li CY, Wang RH, Zhu CG. An approach for designing a developable surface through a given line of curvature. Comput Aided Des 2013; 45:621-7.

[39] Nassar H. A, Rashad A. B, Hamdoon F.M. Ruled surfaces with timelike rulings. Appl. Math. Comput. 147 (2004) 241–253.